\input amstex \documentstyle{amsppt}\magnification=1200 \nologo
\nopagenumbers \NoBlackBoxes  \hsize=16 truecm \vsize = 9 truein
\voffset = 7.5 truemm \hoffset = 0.1 truemm  \topmatter

\title
Hierarchical structure of the family of curves with maximal genus
verifying flag conditions.
\endtitle

\author Vincenzo Di Gennaro \endauthor

\leftheadtext{Vincenzo Di Gennaro} \rightheadtext{Hierarchical
structure of the family of curves with maximal genus}

\medskip
\abstract \nofrills

ABSTRACT. Fix integers $r,s_1,\dots,s_l$ such that $1\leq l\leq
r-1$ and $s_l\geq r-l+1$, and let $\Cal C(r;s_1,\dots,s_l)$ be the
set of all integral, projective and nondegenerate curves $C$ of
degree $s_1$ in the projective space $\bold P^r$, such that, for
all $i=2,\dots,l$, $C$ does not lie on any integral, projective
and nondegenerate variety of dimension $i$ and degree $<s_i$. We
say that a curve $C$ satisfies the {\it{flag condition}}
$(r;s_1,\dots,s_l)$ if $C$ belongs to $\Cal C(r;s_1,\dots,s_l)$.
Define $ G(r;s_1,\dots,s_l)=max\left\{p_a(C):\,C\in \Cal
C(r;s_1,\dots,s_l)\right \}, $ where $p_a(C)$ denotes the
arithmetic genus of $C$. In the present paper, under the
hypothesis $s_1>>\dots>>s_l$, we prove that a curve $C$ satisfying
the flag condition $(r;s_1,\dots,s_l)$ and of maximal arithmetic
genus $p_a(C)=G(r;s_1,\dots,s_l)$ must lie on a unique flag such
as $C=V_{s_1}^{1}\subset V_{s_2}^{2}\subset \dots \subset
V_{s_l}^{l}\subset {\bold P^r}$, where, for any $i=1,\dots,l$,
$V_{s_i}^i$ denotes an integral projective subvariety of ${\bold
P^r}$ of degree $s_i$ and dimension $i$, such that its general
linear curve section satisfies the flag condition
$(r-i+1;s_i,\dots,s_l)$ and has maximal arithmetic genus
$G(r-i+1;s_i,\dots,s_l)$. This proves the existence of a sort of a
hierarchical structure of the family of curves with maximal genus
verifying flag conditions.



\endabstract
\endtopmatter

\bigskip

Fix integers $r,d,s$ such that $s\geq r-1$, and let $\Cal
C(r;d,s)$ be the set of all integral, projective and nondegenerate
curves  of degree $d$ in the projective space $\bold P^r$, not
contained in  any integral, projective surface of degree $<s$.
Extending classical results of Halphen [H], Noether [N] and
Castelnuovo [C], and more recent results of Gruson and Peskine
[GP], and Eisenbud and Harris [EH], in [CCD] one proves that, when
$d>>s$, the curves of maximal arithmetic genus in $\Cal C(r;d,s)$
are contained in surfaces of degree $s$, whose general hyperplane
sections are themselves curves of maximal arithmetic genus in
$\Cal C(r-1;s,r-2)$ (the so called \lq\lq Castelnuovo curves\rq\rq
). In the present paper we show that this property is a particular
case of a more general property.

In order to state our main result, we need some preliminary
notation. Fix integers $r,s_1,\dots,s_l$ such that $1\leq l\leq
r-1$ and $s_l\geq r-l+1$, and let $\Cal C(r;s_1,\dots,s_l)$ be the
set of all integral, projective and nondegenerate curves $C$ of
degree $s_1$ in the projective space $\bold P^r$, such that, for
all $i=2,\dots,l$, $C$ does not lie on any integral, projective
and nondegenerate variety of dimension $i$ and degree $<s_i$. We
say that a curve $C$ satisfies the {\it{flag condition}}
$(r;s_1,\dots,s_l)$ if $C$ belongs to $\Cal C(r;s_1,\dots,s_l)$.
Notice that $\Cal C(r;s_1,r-1)$ is simply the set of all integral,
projective and nondegenerate curves of degree $s_1$ in $\bold P^r$
(i.e. $\Cal C(r;s_1,r-1)=\Cal C(r;s_1)$). Therefore, when studying
the set $\Cal C(r;s_1,\dots,s_l)$, one may assume $l\geq 2$.
Define
$$
G(r;s_1,\dots,s_l)=max\left\{p_a(C):\,C\in \Cal
C(r;s_1,\dots,s_l)\right \},
$$
where $p_a(C)$ denotes the arithmetic genus of $C$. We refer to
[CCD2] for a general discussion on the genus of curves verifying
flag conditions, and its relationship with Castelnuovo-Halphen
Theory. Improving Theorem 3.3 and Corollary 3.4 in [CCD2], in the
present paper we prove the following:
\smallskip
\proclaim{Theorem} Assume that $s_1>>\dots>>s_l$, and fix a curve
$C\in \Cal C(r;s_1,\dots,s_l)$ of maximal arithmetic genus
$p_a(C)=G(r;s_1,\dots,s_l)$. Then one has:

\smallskip
(a) $C$ is arithmetically Cohen-Macaulay;

\smallskip
(b) there exists a unique flag
$$
C=V_{s_1}^{1}\subset V_{s_2}^{2}\subset \dots \subset
V_{s_l}^{l}\subset {\bold P^r},\tag 1.1
$$
where $V_s^j$ denotes an integral projective subvariety of ${\bold
P^r}$ of degree $s$ and dimension $j$;

\smallskip
(c) for any $i=1,\dots,l$ one has
$$
V_{s_i}^{(1)}\in \Cal C(r-i+1;s_i,\dots,s_l)
\quad{\text{and}}\quad p_a(V_{s_i}^{(1)})=G(r-i+1;s_i,\dots,s_l),
$$
where $V_{s_i}^{(1)}$ denotes the curve intersection of
$V_{s_i}^{i}$ with a general linear subspace of ${\bold P^r}$ of
dimension $r-i+1$;

\smallskip
(d) there exists a rational number $R=R(r;s_1,\dots,s_l)$
depending only on $r,s_1,\dots,s_l$, such that
$$
G(r;s_1,\dots,s_l)=
{\frac{s_1^2}{2s_2}}+{\frac{s_1}{2s_2}}\left[2G(r-1;s_2,\dots,s_l)-2-s_2\right]+R
$$
and $|R|\leq s_2^3/(r-2)$.
\endproclaim
\smallskip

Properties (b) and (c) above show a sort of a hierarchical
structure of the family of curves with maximal genus verifying
flag conditions. Moreover, with the exception of the \lq\lq
constant term\rq \rq $R$, property (d) gives a recurrence formula
for $G(r;s_1,\dots,s_l)$. In Remark (iii) below,  we make explicit
what the condition $s_1>>\dots>>s_l$ means.

We will prove  Theorem using some of the results contained in
[CCD2], and using induction on $l$, the case $l=2$ being contained
in the main result of [CCD]. The induction argument relies on the
following:
\smallskip
\proclaim{Lemma} Let $S\subset \bold P^r$ be an irreducible,
reduced, nondegenerate projective surface, of degree $s\geq
r-1\geq 2$. Denote by $S^{(1)}\subset \bold P^{r-1}$ the general
hyperplane section of $S$, by $\pi$ its arithmetic genus, by $\Cal
I_{S^{(1)}}$ its ideal sheaf in $\bold P^{r-1}$, by $S^{(0)}$ the
general hyperplane section of $S^{(1)}$ and by $h_{S^{(0)}}$ its
Hilbert function. For any integer $i$, denote by $\delta_i$ the
dimension of the kernel of the natural map $H^1(\bold P^{r-1},
\Cal I_{S^{(1)}}(i-1))\to H^1(\bold P^{r-1}, \Cal
I_{S^{(1)}}(i))$.

Let $C\subset S$ be an irreducible, reduced, nondegenerate
projective curve, of degree $d\geq s^2+s(r-4)^2$ for $3\leq r\leq
4$, and of degree $d> s^2-s$ for $r\geq 5$. Denote by $p_a(C)$ the
arithmetic genus of $C$, by $\Gamma$ the general hyperplane
section of $C$, by $h_{\Gamma}$ its Hilbert function, and define
$m$ and $\epsilon$ by dividing
$$
d-1=ms+\epsilon,\quad 0\leq \epsilon\leq s-1.
$$
Moreover define:
$$
R(C)={\frac{1+\epsilon}{2s}}(s+1-\epsilon -2\pi)
-\sum_{i=1}^{+\infty}(i-1)(s-h_{S^{(0)}}(i))+
\sum_{i=1}^{+\infty}(i-1)\delta_i+\sum_{i=m+1}^{+\infty}(d-h_{\Gamma}(i)).
$$
Then one has:
$$
\quad  p_a(C)\leq \sum_{i=1}^{+\infty}(d-h_{\Gamma}(i))=
{\frac{d^2}{2s}}+{\frac{d}{2s}}(2\pi-2-s)+R(C) \tag 2.1
$$
and $|R(C)|\leq s^3/(r-2)$.
\endproclaim \bigskip

The proof of Lemma entirely relies on Castelnuovo Theory. In
particular we use the following general formula
$$
\sum_{i=1}^{+\infty}(s-h_{S^{(0)}}(i))=\pi
+\sum_{i=1}^{+\infty}\delta_i
$$
(see [Ci], pg. 31) which enables us to compute the coefficient of
the linear term $d$ in (2.1).

Notice that, when the surface $S\subset {\bold P^r}$ is smooth and
subcanonical, using Hodge Index Theorem and the adjunction
formula, one  has
$$
p_a(C)\leq {\frac{d^2}{2s}}+{\frac{d}{2s}}(2\pi-2-s)+1,
$$
for any curve $C\subset S$. Therefore, one may interpret (2.1) as
a \lq \lq coarse numerical adjunction formula\rq \rq, which holds
for integral projective curves on {\it{any }} integral projective
surface $S\subset {\bold P^r}$.

On the other hand, by [CCD], Main Theorem and Proposition 4.2, we
know that on any such a surface (when $d>s^2-s$) one has
$$
p_a(C)\leq {\frac{d^2}{2s}}+{\frac{d}{2s}}(2G(r-1;s,r-2)-2-s)+R_1,
\tag 0.1
$$
where $G(r-1;s,r-2)$ is the Castelnuovo bound for a nondegenerate
curve of degree $s$ in ${\bold P^{r-1}}$, and $R_1$ is a rational
number which depends only on $s$, $r$ and $\epsilon$ (for the
exact definition of $R_1$ we refer to [CCD], pg. 230-231). On
(certain) Castelnuovo surfaces, i.e. surfaces whose general
hyperplane section is a Castelnuovo curve, previous bound (0.1) is
sharp (see [CCD], pg. 243-244). Now, in view of our Lemma, we may
refine the bound (0.1) proved  in [CCD], in the following sense:

\smallskip
\proclaim{Corollary} Fix integers $r$, $s$ and $\pi$ such that
$r\geq 3$, $s\geq r-1$ and $\pi\geq 0$. Let $C\subset {\bold P^r}$
be an integral, nondegenerate,  projective curve of degree $d>>s$.
Assume that $C$ is not contained in any surface of degree $<s$,
and not contained in any surface of degree $s$ with linear genus
$>\pi$. Then one has:
$$
p_a(C)\leq
{\frac{d^2}{2s}}+{\frac{d}{2s}}(2\pi-2-s)+{\frac{s^3}{r-2}}.\tag
3.1
$$
\endproclaim \bigskip

In Remark (iii) below,  we make explicit  the condition $d>>s$.

Notice that the bound (3.1) given in Corollary is not sharp.
However, dividing (3.1) by $d$, and assuming as before $d>>s$, we
get
$$
e(C)\leq {\frac{d}{s}}+{\frac{2\pi-2-s}{s}}, \tag 0.2
$$
where $e(C)$ denotes the speciality index of $C$, i.e.
$$
e(C)=max\left\{t\in\bold Z :\,H^1(C,\Cal O_C(t))\neq 0\right \}
$$
(recall that by [GP], pg. 51, Remarque 3.6, one has $e(C)d\leq
2p_a(C)-2$). Now, at least in certain cases, previous bound (0.2)
is sharp (e.g. when $C$ is a complete intersection on a complete
intersection surface of degree $s$ and linear genus $\pi$). The
bound (0.2) should be compared with the  \lq\lq Th\'eor\`eme de
sp\'ecialit\'e\rq\rq in [GP], pg. 32. We have in mind to give more
information on  (0.2) in a forthcoming paper.

Now we are going to prove the announced results. We work over the
complex field and we use standard notation of Algebraic Geometry.
We begin by showing the Lemma.

\demo {Proof of Lemma} In view of [GP] we may assume $r\geq 4$.
For the proof of the inequality $p_a(C)\leq
\sum_{i=1}^{+\infty}(d-h_{\Gamma}(i))$ we refer to [EH], Corollary
3.2.

In order to compute the sum
$\sum_{i=1}^{+\infty}(d-h_{\Gamma}(i))$, first notice that by
Bezout's Theorem we have $h_{\Gamma}(i)=h_{S^{(1)}}(i)$ for any
$i\leq m$, where $h_{S^{(1)}}$ denotes the Hilbert function of
$S^{(1)}$. Hence we may write:
$$
\sum_{i=1}^{+\infty}(d-h_{\Gamma}(i))=\sum_{i=1}^{m}(d-h_{S^{(1)}}(i))
+\sum_{i=m+1}^{+\infty}(d-h_{\Gamma}(i)).\tag 2.2
$$
From [Ci], pg.30, we know that
$$
(h_{S^{(1)}}(j)-h_{S^{(1)}}(j-1))-h_{S^{(0)}}(j)=\delta_j
$$
for any integer $j$. Therefore we have:
$$\multline \\
\sum_{i=1}^{m}(d-h_{S^{(1)}}(i))= \sum_{i=1}^{m}\left[
d-\sum_{j=0}^{i}(h_{S^{(1)}}(j)-h_{S^{(1)}}(j-1))\right]=\\
\sum_{i=1}^{m}\left[
d-\sum_{j=0}^{i}(h_{S^{(0)}}(j)+\delta_j)\right]=
md-\left[m(1+\delta_0)+\sum_{i=1}^{m}(m-i+1)(h_{S^{(0)}}(i)+\delta_i)\right]=\\
m(d-1)-\sum_{i=1}^{m}(m-i+1)(s-s+h_{S^{(0)}}(i)+\delta_i)=\\
{\frac{ms}{2}}(m-1)+m\epsilon+m\sum_{i=1}^{m}(s-h_{S^{(0)}}(i)-\delta_i)
-\sum_{i=1}^{m}(i-1)(s-h_{S^{(0)}}(i)-\delta_i).
\endmultline
$$
Now define $w$ and $v$ by dividing
$$
s-1=w(r-2)+v, \quad 0\leq v \leq r-3. \tag 2.3
$$
By [EH], Theorem 3.7, and [GLP], we know that
$$
h_{S^{(0)}}(i)=s\quad {\text{for any $i\geq w+1$}},
$$
and
$$
H^1(\bold P^{r-1}, \Cal I_{S^{(1)}}(i))=0\quad {\text{for any
$i\geq s-r+2$}}. \tag 2.4
$$
Since $d>>s$ then $m\geq w+1$ and $m\geq s-r+2$. And so, taking
into account [Ci], pg. 31, we get
$$
\sum_{i=1}^{m}(s-h_{S^{(0)}}(i)-\delta_i)=
\sum_{i=1}^{+\infty}(s-h_{S^{(0)}}(i)-\delta_i)=\pi
$$
and
$$
\sum_{i=1}^{m}(i-1)(s-h_{S^{(0)}}(i)-\delta_i)=\sum_{i=1}^{+\infty}(i-1)(s-h_{S^{(0)}}(i)-\delta_i).
$$
Continuing previous computation, we have
$$
\sum_{i=1}^{m}(d-h_{S^{(1)}}(i))=
{\frac{ms}{2}}(m-1)+m\epsilon+m\pi
-\sum_{i=1}^{+\infty}(i-1)(s-h_{S^{(0)}}(i)-\delta_i).
$$
Replacing $m$ with $(d-1-\epsilon)/s$, and taking into account
(2.2), we get (2.1).

To conclude the proof of Lemma  we have to estimate $R(C)$. We
will analyze each of the four terms appearing in the definition of
$R(C)$.

First notice that
$$
-\pi\leq {\frac{1+\epsilon}{2s}}(s+1-\epsilon -2\pi)\leq
{\frac{s+1}{2}}.
$$
We may estimate the arithmetic genus $\pi$ of $S^{(1)}$ using
Castelnuovo bound for curves of degree $s$ in $\bold P^{r-1}$
(compare with (2.3)):
$$
\pi\leq {{w}\choose{2}}(r-2)+wv\leq {\frac{s^2}{2(r-2)}},
$$
from which we obtain
$$
-{\frac{s^2}{2(r-2)}}\leq {\frac{1+\epsilon}{2s}}(s+1-\epsilon
-2\pi)\leq {\frac{s+1}{2}}. \tag 2.5
$$

Now we turn to next term. From [EH] we know that
$$
h_{S^{(0)}}(i)\geq min\{s,\, i(r-2)+1\} \tag 2.6
$$
for any $i\geq 1$. Since
$$
\sum_{i=1}^{+\infty}(i-1)(s-min\{s,\, i(r-2)+1\})=
{{w}\choose{3}}(r-2)+{{w}\choose{2}}v\leq {\frac{s^3}{3(r-2)^2}},
$$
it follows that
$$
0\leq \sum_{i=1}^{+\infty}(i-1)(s-h_{S^{(0)}}(i))\leq
{\frac{s^3}{3(r-2)^2}}. \tag 2.7
$$

In order to estimate the third term, first notice that from (2.4)
we have
$$
\sum_{i=1}^{+\infty}(i-1)\delta_i= \sum_{i=1}^{s}(i-1)\delta_i\leq
(s-1)\sum_{i=1}^{+\infty}\delta_i.
$$
On the other hand, from [Ci], pg. 31, and (2.6), we have
$$
\sum_{i=1}^{+\infty}\delta_i\leq
\sum_{i=1}^{+\infty}(s-h_{S^{(0)}}(i))\leq
\sum_{i=1}^{+\infty}(s-min\{s,\, i(r-2)+1\})=
{{w}\choose{2}}(r-2)+wv\leq {\frac{s^2}{2(r-2)}}.
$$
Putting all together we get
$$
0\leq \sum_{i=1}^{+\infty}(i-1)\delta_i\leq
{\frac{s^2(s-1)}{2(r-2)}}. \tag 2.8
$$

Finally we are going to analyze the last term. From [CCD],
Propositions 4.1 and 4.2, we know that
$$
\sum_{i=m+1}^{+\infty}(d-h_{\Gamma}(i))=
\sum_{i=m+1}^{m+w}(d-h_{\Gamma}(i)).
$$
Since $h_{\Gamma}(i)\geq h_{\Gamma}(m)$ for $i\geq m$, we deduce
that
$$
\sum_{i=m+1}^{+\infty}(d-h_{\Gamma}(i))\leq w(d-h_{\Gamma}(m)).
\tag 2.9
$$
Since $d>>s$, by Bezout's Theorem and (2.4) we have
$$
h_{\Gamma}(m)=h_{S^{(1)}}(m))=h^0(S^{(1)}, \Cal O_{S^{(1)}}(m)).
\tag 2.10
$$
On the other hand, from (2.6) we deduce
$$
h^1({\bold P^{r-2}}, \Cal I_{S^{(0)}}(i))=0
$$
for any $i\geq w+1$ ($\Cal I_{S^{(0)}}$=ideal sheaf of $S^{(0)}$
in ${\bold P^{r-2}}$). This implies that
$$
h^1(S^{(1)}, \Cal O_{S^{(1)}}(m))=0.
$$
From (2.10) it follows that
$$
h_{\Gamma}(m)=ms+1-\pi.
$$
From  Castelnuovo's bound on $\pi$ we deduce
$$
w(d-h_{\Gamma}(m))=w(\epsilon +\pi)\leq {\frac{s^3}{2(r-2)^2}},
$$
and so, from (2.9), we get
$$
0\leq \sum_{i=m+1}^{+\infty}(d-h_{\Gamma}(i))\leq
{\frac{s^3}{2(r-2)^2}}. \tag 2.11
$$

Using (2.5), (2.7), (2.8) and (2.11), we obtain the estimate for
$|R(C)|$. This concludes the proof of Lemma. $\square$
\bigskip
\bigskip
Next we give the proof of Theorem.

\demo {Proof of Theorem} The case $l= 2$ is contained in the main
result of [CCD]. Therefore we may argue by induction on $l$, and
assume $l\geq 3$.

For the existence of the flag (1.1) we refer to [CCD2], Corollary
2.8. The uniqueness follows by Bezout's Theorem and the assumption
$s_1>>\dots>>s_l$. This proves property (b).

By Lemma  we know that
$$
p_a(C)\leq
{\frac{s_1^2}{2s_2}}+{\frac{s_1}{2s_2}}\left(2\pi-2-s_2\right)+s_2^3/(r-2),\tag
1.2
$$
where $\pi$ denotes the linear arithmetic genus of $V_{s_2}^{2}$.
Since $s_2>>\dots>>s_l$, then by Bezout's Theorem we have
$V_{s_2}^{(1)}\in \Cal C(r-1;s_2,\dots,s_l)$, and so
$$
\pi \leq G(r-1;s_2,\dots,s_l). \tag 1.3
$$

Now fix a curve $D\in \Cal C(r-1;s_2,\dots,s_l)$ of maximal
arithmetic genus. By induction, this curve is arithmetically
Cohen-Macaulay and determines a flag
$$
D=W_{s_2}^{1}\subset W_{s_3}^{2}\subset \dots \subset
W_{s_l}^{l-1}\subset {\bold P^{r-1}}.\tag 1.4
$$
Since $s_1>>s_2$, by [CCD2], Lemma 2.6, we may construct on the
cone $C(D)$ over $D$ in ${\bold P^{r}}$, an integral,
nondegenerate, projective and arithmetically Cohen-Macaulay curve
$E\subset C(D)$ of degree $s_1$. $E$ lies on the cone of the flag
(1.4), therefore $E\in \Cal C(r;s_1,s_2,\dots,s_l)$. Moreover,
since the general hyperplane section of $C(D)$ has arithmetic
genus $G(r-1;s_2,\dots,s_l)$ and $E$ is arithmetically
Cohen-Macaulay, then by [EH], Remark 3.1.1, and our Lemma, we have
$$
p_a(E)=\sum_{i=1}^{+\infty}(s_1-h_{E'}(i))=
{\frac{s_1^2}{2s_2}}+{\frac{s_1}{2s_2}}\left[2G(r-1;s_2,\dots,s_l)-2-s_2\right]+R_2,\tag
1.5
$$
with $|R_2|\leq s_2^3/(r-2)$ ($h_{E'}$= Hilbert function of the
general hyperplane section $E'$ of $E$). Since $p_a(E)\leq
G(r;s_1,\dots,s_l)=p_a(C)$ and $|R_2|\leq s_2^3/(r-2)$, and since
$s_1>>s_2$, then from (1.2), (1.3) and (1.5) we get
$$
\pi=G(r-1;s_2,\dots,s_l).\tag 1.6
$$
This means that $V_{s_2}^{(1)}$ is a curve of maximal genus
verifying the flag condition $(r-1;s_2,\dots,s_l)$. This proves,
by induction, property (c).

In particular, $V_{s_2}^{(1)}$ is arithmetically Cohen-Macaulay.
Hence, using again [CCD2], Lemma 2.6, we may construct an
arithmetically Cohen-Macaulay curve $F$ belonging to $\Cal
C(r;s_1,\dots,s_l)$, whose general hyperplane section $F'$ has the
same Hilbert function as the general hyperplane section $C'$ of
$C$. It follows that
$$\split
G(r;s_1,\dots,s_l)=p_a(C)\leq
\sum_{i=1}^{+\infty}(s_1-h_{C'}(i))\\
= \sum_{i=1}^{+\infty}(s_1-h_{F'}(i))=p_a(F')\leq
G(r;s_1,\dots,s_l),
\endsplit
$$
and so
$$
p_a(C)=\sum_{i=1}^{+\infty}(s_1-h_{C'}(i)),\tag 1.7
$$
i.e. $C$ is arithmetically Cohen-Macaulay. This proves property
(a).

At this point, taking into account that $C\subset V_{s_2}^{2}$,
property (d) follows from (1.7) and previous Lemma. This concludes
the proof of Theorem. $\square$
\bigskip
\bigskip

Finally we turn to the proof of Corollary.

\demo{Proof of Corollary} If $C$ is not contained in any surface
of degree $<s+1$, then, from the main result of [CCD] and our
Lemma (compare with (0.1)), we deduce that
$$
p_a(C)\leq
{\frac{d^2}{2(s+1)}}+{\frac{d}{2(s+1)}}(2G(r-1;s+1,r-2)-2-(s+1))+{\frac{(s+1)^3}{r-2}},\tag
3.2
$$
where $G(r-1;s+1,r-2)$ is the Castelnuovo bound for a
nondegenerate curve of degree $s+1$ in ${\bold P^{r-1}}$. Since
$d>>s$, previous bound (3.2) is strictly less than the bound
appearing in (3.1). Therefore we may assume that $C$ is contained
on some surface of degree $s$, with linear genus $\leq \pi$. In
this case Corollary follows from Lemma. This concludes the proof
of Corollary. $\square$
\bigskip \bigskip

\noindent {\bf Remark.} (i) With the same notation as in Lemma, we
notice that when the surface $S$ is arithmetically Cohen-Macaulay,
then all $\delta_i$ vanish and
$\sum_{i=1}^{+\infty}(i-1)(s-h_{S^{(0)}}(i))$ is equal to the
arithmetic genus $p_a(S)$ of $S$ (see [D], Remark 2.3). Therefore,
in this case, we have
$$
R(C)={\frac{1+\epsilon}{2s}}(s+1-\epsilon -2\pi) -p_a(S)
+\sum_{i=m+1}^{+\infty}(d-h_{\Gamma}(i)).
$$
In particular, in Theorem, since we know that the surface
$V_{s_2}^2$ is arithmetically Cohen-Macaulay,  we have
$$
R={\frac{1+\epsilon}{2s}}(s+1-\epsilon -2G(r-1;s_2,\dots,s_l))
-p_a(V_{s_2}^2)+\sum_{i=m+1}^{+\infty}(d-h_{\Gamma}(i)),
$$
where $h_{\Gamma}$ is the Hilbert function of the general
hyperplane section of any maximal curve $C\in\Cal
C(r;s_1,\dots,s_l)$.

(ii) Again in Lemma, we notice that when $S$ is a Castelnuovo
surface,  using the main result of [CCD], one may prove  that
$$
R(C)=O(s^2).
$$

(iii) In proving Theorem, we need the numerical assumption
$s_1>>\dots>>s_l$ only  to use Corollary 2.8 in [CCD2],   Bezout's
Theorem, and to prove (1.6). To this purpose, it suffices to
assume, for any $i=1,\dots,l-1$,
$$
s_i\geq 8(l-1)\left [(l-i+1)^2+2(l-i+1)+9\right
]{\frac{(s_{i+1}+1)^3}{r-i-1}},
$$
$$
s_i>{\frac{(s_{i+1}+1)^2}{r-i-1}}+(2r-2){(s_{i+1}+1)},
$$
$$
s_i>2{\frac{(s_{i+1}+1)}{r-i-1}}\prod_{j=1}^{r-1-i}[(r-i)!(s_{i+1}+1)]^{\frac{1}{r-i-j}},
\quad {\text{and}}
$$
$$
s_i>2{\frac{s_{i+1}^4}{r-i-1}}.
$$

We also may explicit the numerical assumption $d>>s$ made in
Corollary. In fact, we only need it for using [CCD], and to
compare (3.2) with (3.1). To this aim, it suffices to assume
$$
d>{\frac{2(s+1)}{r-2}}\prod_{i=1}^{r-2}[(r-1)!(s+1)]^{\frac{1}{r-1-i}}\quad
{\text{and}}\quad d>{\frac{6(s+1)^3}{r-2}}.
$$
\bigskip \bigskip

{
\bigskip

\noindent Author address:\par \noindent Vincenzo Di Gennaro\par
\noindent Universit\`a di Roma Tor Vergata, Dipartimento di
Matematica,\par \noindent Via della Ricerca Scientifica, 00133
Roma, Italia.\par
\medskip
\noindent E-mail: digennar\@axp.mat.uniroma2.it

\end